\documentclass[12pt]{article}

\usepackage{amsmath}
\usepackage{amsfonts}

\newtheorem{theorem}{Theorem}
\newtheorem{lemma}[theorem]{Lemma}
\newtheorem{corollary}[theorem]{Corollary}
\newtheorem{definition}[theorem]{Definition}
\newtheorem{remark}[theorem]{Remark}

\begin{document}

\title{Pseudodifferential $p$-adic vector fields \\ and pseudodifferentiation of a composite $p$-adic function}

\author{S.Albeverio\footnote{Universit\"at Bonn,
Institut f\"ur Angewandte Mathematik, Bonn, Germany},
S.V.Kozyrev\footnote{Steklov Mathematical Institute, Moscow, Russia}}

\maketitle

\begin{abstract}
We discuss transformation of $p$-adic pseudodifferential operators (in the one--dimensional and multidimensional cases) with respect to $p$-adic maps which correspond to automorphisms of the tree of balls in the corresponding $p$-adic spaces.

In the dimension one we find a rule of transformation for pseudodifferential operators. In particular we find the formula of pseudodifferentiation of a composite function with respect to the Vladimirov $p$-adic fractional operator. We describe the frame of wavelets for the group of parabolic automorphisms of the tree ${\cal T}(\mathbb{Q}_p)$ of balls in $\mathbb{Q}_p$.

In many dimensions we introduce the group of ${\rm mod}\, p$-affine transformations, the family of pseudodifferential operators corresponding to  pseudodifferentiation along vector fields on the tree ${\cal T}(\mathbb{Q}_p^d)$ and obtain a rule of transformation of the introduced pseudodifferential operators with respect to ${\rm mod}\, p$-affine transformations.

\end{abstract}


\section{Introduction}

The present paper discusses action of some $p$-adic groups of transformations on $p$-adic pseudodifferential operators.

$p$-Adic pseudodifferential operators (see the Appendix) were discussed in the books \cite{VVZ,Kochubei} and in many papers, see in particular \cite{wavelets,nhoper,Izv,ACHA,MathSbornik,framesdimd}. These operators are non--local. In the real case the analogous operators possess only a finite dimensional group of symmetries. For applications of $p$-adic pseudodifferential operators and other models of $p$-adic mathematical physics see \cite{obzor}, \cite{Andr3}.

$p$-Adic wavelets were introduced in \cite{wavelets} and discussed in many papers, in particular, multidimensional wavelet bases and relation to representations of some $p$-adic groups of transformations were considered in \cite{framesdimd}. $p$-Adic wavelets are associated to balls in $\mathbb{Q}_p$ --- to any ball we put in correspondence the wavelets supported in this ball which are constant on maximal subballs of the ball. The set of balls in $\mathbb{Q}_p$ can be considered as a tree ${\cal T}(\mathbb{Q}_p)$ (see the Appendix for the definitions).

We can consider the linear combinations of wavelets (i.e. the $p$-adic test functions from the space $D_0(\mathbb{Q}_p)$ of mean zero locally constant complex valued functions with compact support) as sections of the fiber bundle over the tree ${\cal T}(\mathbb{Q}_p)$. These sections take values in the finite dimensional wavelet spaces associated to balls in $\mathbb{Q}_p$.

Since $p$-adic wavelets are eigenvectors of $p$-adic pseudodifferential operators \cite{wavelets,nhoper,Izv,ACHA,MathSbornik}  $p$-adic pseudodifferential operators can be considered as {\it local} operators in the wavelet fiber bundle on the tree ${\cal T}(\mathbb{Q}_p)$. In this language a pseudodifferential operator is a vector field over the tree ${\cal T}(\mathbb{Q}_p)$.

Therefore it is natural to investigate transformations of  $p$-adic pseudodifferential operators with respect to automorphisms of the tree ${\cal T}(\mathbb{Q}_p)$. These transformations can be considered as analogs of transformations of vector fields with respect to diffeomorphisms in differential geometry. In this paper we consider some basic objects of the $p$-adic pseudodifferential geometry, related to wavelets, pseudodifferential operators and automorphisms of trees.

We consider the one--dimensional and multidimensional cases.

In the one--dimensional case we investigate the action of the group generated by isometries and dilations. This group can be considered as a group of automorphisms of the tree ${\cal T}(\mathbb{Q}_p)$ of balls in $\mathbb{Q}_p$  which conserve the infinite point of the absolute of the tree (the parabolic subgroup of the group of automorphisms of the tree ${\cal T}(\mathbb{Q}_p)$). Groups of automorphisms of trees and their representations were studied in \cite{Olshansky,Neretin,Serre,Serre1,Cartier1,Cartier2}. The discussed groups do not possess finite sets of generators. In this sense these groups are infinite dimensional.

In the present paper we prove the following statement: let $\phi$ be a one--to--one map $\phi:
\mathbb{Q}_p\to\mathbb{Q}_p$, which satisfies the following property: the image and the inverse image of any ball are balls (i.e. the map is a ball--morphism, see for example \cite{Neretin}).
Let the map $\phi$ also be continuously differentiable.

The corresponding map of functions is defined as follows:
\begin{equation}\label{on_functions}
\Phi f(x)=f(\phi(x)),\quad x\in\mathbb{Q}_p,\quad f:\mathbb{Q}_p\to\mathbb{C}.
\end{equation}

We show that the Vladimirov operator $D^{\alpha}$ of
$p$-adic fractional differentiation will satisfy the following formula of pseudodifferentiation of a composite function
$$
D^{\alpha}\circ \Phi f(x)=|\phi'(x)|_p^{\alpha}\,\Phi\circ
D^{\alpha}f(x).
$$

We discuss the action of this group in the basis of wavelets and show that the orbit of the unitary action of the group of ball--morphisms is a tight frame (see the Appendix for definitions). We compute a bound for this frame.

In the multidimensional case we perform the following constructions. We put in correspondence to subballs of a $d$-dimensional $p$-adic ball the elements of the module $\mathbb{F}_p^{d}$ (the $d$-dimensional module over the field of residues ${\rm mod}\, p$). Taking the union over all balls in  $\mathbb{Q}_p^d$ we obtain a natural fiber bundle over the base ${\cal T}(\mathbb{Q}_p^d)$ with layers $\mathbb{F}_p^{d}$.

We introduce the group of ${\rm mod}\, p$-affine isometries.
Transformations from this group act as affine transformations in the described above fiber bundle (i.e. the transformations are automorphisms of ${\cal T}(\mathbb{Q}_p^d)$ and the corresponding maps of the set of subballs of a ball to the set of subballs of the image of the ball act as affine transformations on the layers $\mathbb{F}_p^{d}$ of the fiber bundle).

We consider the sections of the introduced $\mathbb{F}_p^{d}$-valued fiber bundle over ${\cal T}(\mathbb{Q}_p^d)$ (a $\mathbb{F}_p^{d}$-valued vector fields). To any $\mathbb{F}_p^{d}$-valued vector field we put in correspondence some pseudodifferential operator acting in $L^2(\mathbb{Q}_p^d)$ (a pseudodifferentiation along $\mathbb{F}_p^{d}$-valued vector field on ${\cal T}(\mathbb{Q}_p^d)$). The introduced family of pseudodifferential operators is new.

Then we find a rule of transformation of the introduced pseudodifferential operators with respect to ${\rm mod}\, p$-affine isometries. In particular we show that a pseudodifferential operator from the introduced family maps to an operator from the same family.

The exposition of the present paper is as follows.

In Section 2 we discuss groups of automorphisms of trees related to isometries and dilations of $\mathbb{Q}_p$. These results in general may be found in the literature but we put here the discussion to give the exposition of relation of $p$-adic analysis and groups of automorphisms of trees.

In Section 3 using the results of Section 2 we describe transformations of (one--dimensional) $p$-adic pseudodifferential operators
with respect to maps $\mathbb{Q}_p\to\mathbb{Q}_p$ corresponding to parabolic automorphisms of the tree ${\cal T}(\mathbb{Q}_p^d)$.

In Section 4 we describe the frame in $L^2(\mathbb{Q}_p)$ obtained by the action of the parabolic group on wavelets.

In Section 5 we give the definition of ${\rm mod}\, p$-affine isometries of $\mathbb{Q}_p^d$.

In Section 6 we describe a relation between vectors in the above mentioned fiber bundle over the tree ${\cal T}(\mathbb{Q}_p^d)$ with the layer $\mathbb{F}_p^{d}$ and some families of subset in $\mathbb{Q}_p^d$.

In Section 7 we introduce a family of multidimensional pseudodifferential operators related to vector fields the tree ${\cal T}(\mathbb{Q}_p^d)$
with values in $\mathbb{F}_p^{d}$ and find a formula of transformation of the introduced operators with respect to  ${\rm mod}\, p$-affine isometries f $\mathbb{Q}_p^d$.

In Section 8 (the Appendix) we put some results of $p$-adic analysis.

\section{Groups acting on trees}

There exist a known duality between ultrametric spaces and trees (see the Appendix for discussion). In the $p$-adic case the corresponding tree can be considered as a Bruhat--Tits tree. The boundary of the Bruhat--Tits tree is a $p$-adic projective line (absolute of the tree), and the tree can be considered as a tree of balls in $\mathbb{Q}_p$. We consider the group of automorphisms of the Bruhat--Tits tree (or ball--morphisms).

Let us remind that the absolute of the tree can be considered as a set of equivalence classes of half--infinite paths in the tree, the paths are equivalent if they coincide starting from some vertex of the tree (i.e. an equivalence class contains paths with coinciding tails).

Consider the following subgroups in the group of automorphisms of the Bruhat--Tits tree.

\medskip

1) The parabolic subgroup contains automorphisms of the tree which conserve some fixed point of the absolute. This point can identified with the infinite point $\infty$ of the projective line. Therefore the parabolic group can be identified with the group of ball--morphisms of $\mathbb{Q}_p$, i.e. the group of one--to--one maps of $\mathbb{Q}_p$ for which the image and the inverse image of any ball are balls.

\medskip

2) The orispheric subgroup in the parabolic group contains ball--morphisms $\phi$ which conserve not only the point  $\infty$ of the absolute of the tree, but for any $\phi$ in the orispheric subgroup there exist a path in the tree in the equivalence class of $\infty$ such that $\phi$ conserves the tail of the path (i.e. maps any vertex in the path to itself starting from some vertex in the path). Equivalently, the orispheric map $\phi$ conserves some ball $I$ in $\mathbb{Q}_p$ and the sequence of increasing balls starting in $I$.

\begin{lemma}\label{isometry}\quad{\sl
The orispheric group coincides with the group of isometries in
$\mathbb{Q}_p$. }
\end{lemma}

\noindent{\it Proof}\quad Let us show that the transformation from the orispheric group is an isometry. Let the transformation $\phi$ from the orispheric group conserve some ball $I$ (and therefore also conserve the sequence of increasing balls starting from $I$).

For a ball $J$ we consider the image $\phi(J)$ of this ball. Let us show that
$$
{\rm diam}\,(J)={\rm diam}\,(\phi(J)),
$$
where ${\rm diam}$ is the diameter of the ball.

By definition
$$
\phi({\rm sup}(I,J))={\rm sup}(I,J),
$$
where ${\rm sup}(I,J)$ is the minimal ball which contains the both balls $I$, $J$.

Consider the maximal (finite) sequence of increasing embedded balls $\{J_1,\dots,J_{m}\}$, $J_1=J$, $J_m={\rm
sup}(I,J)$, $J_{i}\subset J_{i+1}$. The diameters of the neighbor (embedded) balls in this sequence
differ by a multiplication by $p$.

The image of this sequence is
$\{\phi(J_1),\dots,\phi(J_{m})\}$, $\phi(J_1)=\phi(J)$, where $\phi(J_m)=\phi({\rm
sup}(I,J))={\rm sup}(I,J)=J_m$ (by orisphericity of the map $\phi$).

Since $\phi$ is a ball--morphism, the sequence $\{\phi(J_1),\dots,\phi(J_{m})\}$ is a maximal sequence of nested balls (of the same length $m$). The diameters of the neighbor balls in this sequence differ by a multiplication by $p$. Since
$$
{\rm diam}\,(J_m)={\rm diam}\,(\phi(J_m)),
$$
we obtain
$$
{\rm diam}\,(J)={\rm diam}\,(\phi(J)),
$$
i.e. $\phi$ is an isometry.

\medskip

Let us show that if a map $\phi$ is an isometry then $\phi$ belongs to the orispheric group.

Let us show that an isometry maps a ball into a ball with the same diameter. Consider a ball
$I$ and $x\in I$ is a point (a center of $I$). Then the image $\phi(I)$
is a subset of a ball $J$, ${\rm diam}\,(J)={\rm diam}\,(I)$, and
$\phi(x)$ is a center of $J$.

By isometricity of $\phi$ the set $\phi(I)$ contains points which belong to at least two different maximal subballs in $J$. The inverse images of these points (the inverse images are unique by isometricity) will belong to the different maximal subballs in  $I$. If $I$ contains other maximal subball, we choose a point in this subball and consider the image of this point with respect to $\phi$. By isometricity of $\phi$ this image will belong to the third (different from the considered earlier) maximal subballs in $J$.

Repeating this argument and taking into account that the ball $J$ contains a finite number of maximal subballs we prove that $\phi(I)$ contains points of all maximal subballs in $J$.

Repeating the above arguments for maximal subballs in maximal subballs $I$ and so forth we get that $\phi(I)$ contains dense subset of $J$. Then by continuity of the isometry $\phi$ we get that $\phi(I)=J$ (the image w.r.t. $\phi$ of a closed set is closed).

Let us show that the isometry $\phi$ is surjective (and thus is a one--to one map). Assume that there exists $x\in \mathbb{Q}_p$ which does not belong to the image of $\phi$. Then there exists a sufficiently large ball $J$ with a center in
$x$ which contains a point $y$ belonging to the image of $\phi$.

Consider the inverse image $\phi^{-1}(y)$ and a ball $I$ with the diameter ${\rm diam}\,(J)$
and the center $\phi^{-1}(y)$. Then we can apply the above arguments to the pair of balls $I$, $J$ and show that $\phi(I)=J$ (therefore $x$ belongs to the image of $\phi$).

We have proved that an isometry $\phi$ is a ball--morphism which conserves the diameters of balls.

Let us show that for an isometry $\phi$ there exists a ball $I$ with the following property: the map $\phi$ conserves the increasing sequence of balls which begins in $I$ (i.e. the ball--morphism $\phi$ will be orispheric).

Namely for a point $x$ we consider the ball $I={\rm
sup}\,(x,\phi(x))$ (i.e. the minimal ball which contains $x$ and
$\phi(x)$).

Since $I$ contains $x$ and $\phi(x)$, the image $T(I)$ does intersect the ball
$I$. Since ${\rm diam}(I)={\rm diam}(\phi(I))$ and $\phi$ is a ball--morphism we have $I=\phi(I)$. By the same argument for any ball $J\supset I$ we get $\phi(J)=J$.

Therefore the isometry $\phi$ is an orispheric ball--morphism which finishes the proof of the lemma. $\square$

\bigskip

The following statement can be found in \cite{Olshansky}.

\begin{lemma}\quad{\sl
The orispheric group is a normal subgroup in the parabolic group, moreover the factor--group is isomorphic to $\mathbb{Z}$. }
\end{lemma}

For the $p$-adic case, the mentioned above factor--group can be identified with the group of multiplications by degrees of $p$. In particular, we will prove the following lemma.

\begin{lemma}\label{product}\quad{\sl An arbitrary transformation belonging to the parabolic group can be uniquely expressed as a product of  the multiplication by a degree of $p$ and an isometry (i.e. a map from the orispheric group). }
\end{lemma}

\noindent{\it Proof}\quad Let a ball--morphism $\eta$ (i.e. a map from the parabolic group) map some ball $I$ to the ball $\eta(I)$ with the diameter
$p^{\gamma}\,{\rm diam}\,(I)$.

It is easy to see that the multiplication by $p^{\gamma}$ belongs to the parabolic group. Let us show that the (parabolic) ball--morphism $\phi=p^{\gamma}S$ belongs to the orispheric group.

By definition ${\rm diam}\,(\phi(I))={\rm diam}\,(I)$.

Let us consider the ball $J={\rm sup}(I,\phi(I))$ (the minimal ball which contains $I$ and $\phi(I)$). Since ${\rm diam}\,(\phi(I))={\rm diam}\,(I)$ then the lengths of sequences of nested balls between $J$ and $I$, respectively $J$ and $\phi(I)$ are equal.

Since $J={\rm sup}(I,\phi(I))$ then $\phi(J)={\rm sup}(\phi(I),\phi^2(I))$.
Since $\phi$ is a ball--morphism the lengths of sequences of nested balls between $J$ and $I$, $\phi(I)$, and between $\phi(J)$ and $\phi(I)$, $\phi^2(I)$ are equal. Since ${\rm diam}\,(I)={\rm diam}\,(\phi(I))$ we obtain
${\rm diam}\,(J)={\rm diam}\,(\phi(J))$.

By the construction, the ball $\phi(J)$ has non--zero intersection with the ball $J$, i.e. one of the balls $J$, $\phi(J)$ contains the other ball. Since the mentioned balls have equal diameters we get $J=\phi(J)$. Analogously $L=\phi(L)$ for any
$L\supset J$.

We have proved that the map $\phi=p^{\gamma}\eta$ belongs to the orispheric group.

The uniqueness of the obtained decomposition of the parabolic map $S$ in the composition of the orispheric map $\phi$ and dilation by $p^{-\gamma}$ is obvious. $\square$

\bigskip

In contrast to the real case, in the $p$-adic case the group of isometries does not possess a finite set of generators.
The parabolic group (i.e. the group of ball morphisms) for $p$-adic analysis should play the role of the group of diffeomorphisms in real analysis.

One can consider $p$-adic diffeomorphisms.
It is easy to check that there exist diffeomorphisms which do not belong to the group of ball--morphisms.
Also there exist non--differentiable ball--morphisms.

The following lemma can be found in \cite{Schikhof}.

\begin{lemma}\quad{\sl
1) If the function $f:\mathbb{Q}_p\to\mathbb{Q}_p$ is differentiable at $a$ and the derivative at $a$ is not equal to zero, then there exists a ball with a sufficiently small diameter and the center in $a$, satisfying the property: for any $x$ in this ball
$$
|f(x)-f(a)|_p=|f'(a)|_p|x-a|_p.
$$

2) If the function $f:\mathbb{Q}_p\to\mathbb{Q}_p$ is differentiable at $a$ and the derivative at $a$ is equal to zero, then $\forall C>0$
there exists a ball with a sufficiently small diameter and the center in $a$, satisfying the property: for any $x$ in this ball
$$
|f(x)-f(a)|_p\le C|x-a|_p.
$$

}
\end{lemma}

\begin{corollary}\label{derivative_one}\quad{\sl
For a continuously differentiable isometry on $\mathbb{Q}_p$  the norm of the derivative is equal to one. }
\end{corollary}

\section{Pseudodifferentiation of a composite function }

In the present Section we investigate transformations of pseudodifferential operators with respect to the action of parabolic maps.

Consider a pseudodifferential operator of the investigated in \cite{nhoper}, \cite{Izv}, \cite{ACHA}, \cite{MathSbornik}  form
$$
D_{F(\cdot)}f(x)=\int F({\rm sup}(x,y))(f(x)-f(y))dy.
$$
Here the integration kernel $F(I)$ is the complex valued function of balls $I\in{\cal T}(\mathbb{Q}_p)$.

Apply the parabolic transformation $\phi$ which corresponds to the dilation by $p^{\gamma}$:
$$
\Phi\circ D_{F(\cdot)}f(x)=Df(\phi(x))=\int F({\rm sup}(\phi(x),y))(f(\phi(x))-f(y))dy=
$$
$$
=\int F(\phi({\rm sup}(x,y)))(f(\phi(x))-f(\phi(y)))d\phi(y)=p^{-\gamma}\int F(\phi({\rm sup}(x,y)))(f(\phi(x))-f(\phi(y)))dy.
$$

Here $\Phi$ is the map of functions (\ref{on_functions}) which corresponds to the isometry $\phi$.
Since
$$
D_{F(\cdot)}\circ \Phi f(x)=\int F({\rm sup}(x,y))(f(\phi(x))-f(\phi(y)))dy
$$
we get the transformation rule
\begin{equation}\label{transformD_F}
\Phi\circ D_{F(\cdot)}=p^{-\gamma}D_{F(\phi(\cdot))}\circ \Phi.
\end{equation}

Applying this formula to the Vladimirov operator
$$
D^{\alpha}f(x)=\Gamma_p^{-1}(-\alpha)\int
{f(x)-f(y)\over
|x-y|_p^{1+\alpha}}dy
$$
of $p$-adic fractional differentiation we get
a natural formula for the pseudodifferentiation of a composite function:
$$
D^{\alpha}\circ \Phi f(x)=p^{-\gamma\alpha}\,\Phi\circ
D^{\alpha}f(x).
$$

In the case where the parabolic ball--morphism $\phi$ is continuously differentiable on
$\mathbb{Q}_p$, by corollary  \ref{derivative_one} the multiplier
$p^{-\gamma\alpha}$ can be expressed as the degree of the norm of the derivative of the map $\phi$:
$$
D^{\alpha}\circ \Phi f(x)=|\phi'(x)|_p^{\alpha}\,\Phi\circ
D^{\alpha}f(x).
$$

The above formula can be compared with the formula of differentiation of a composite function in real analysis
$$
{d\over dx}f(\phi(x))={df(y)\over
dy}\biggr|_{y=\phi(x)}{d\phi(x)\over dx}.
$$

We have proved that in  $p$-adic analysis for the non-local Vladimirov operator $D^{\alpha}$ of fractional differentiation there exists an analog of the formula of differentiation of a composite function. Let us note that the group of ball--morphisms does not possess a finite set of generators. In this sense this group is infinite dimensional. From the point of view of the analysis of $p$-adic pseudodifferential operators the group of ball--morphisms is a natural analog of the group of diffeomorphisms in real analysis.

\section{Action of parabolic group on wavelets}

The action of the group of (parabolic) ball--morphisms in $\mathbb{Q}_p$ generates a representation of this group in $L^2(\mathbb{Q}_p)$. We consider the unitary representation of the group, i.e. for a ball--morphism equal to a composition of an isometry $\phi$ and a dilation by $p^{\gamma}$ the corresponding transformation in $L^2(\mathbb{Q}_p)$ will have the form
$$
f(x)\mapsto p^{{\gamma\over 2}}f\left(\phi(p^{\gamma}x)\right).
$$

Let us remind that for any ball in $\mathbb{Q}_p$ we have the $p-1$-dimensional space of wavelets with the support equal to this ball.
In particular the wavelet $\psi=\chi(p^{-1}x)\Omega(|x|_p)$ is related to the unit ball  $\mathbb{Z}_p$.

A ball--morphism $\phi$ in $\mathbb{Q}_p$ maps a ball to a ball and maximal subballs in a ball to maximal subballs in a ball. The condition of zero mean is conserved by a ball--morphism. Therefore a wavelet related to a ball $I$ maps to a function related to a ball $\phi(I)$, where the space of functions related to a ball $B$ is defined as the set of mean zero linear combinations of characteristic functions of maximal subballs in $B$.

We have the following lemma about the orbit (see the Appendix for the definition of a frame):

\begin{lemma}\quad{\sl
The orbit of the $p$-adic wavelet $\psi=\chi(p^{-1}x)\Omega(|x|_p)$
with respect to the group of (parabolic) ball--morphisms in $\mathbb{Q}_p$ is the tight uniform frame in $L^2(\mathbb{Q}_p)$ which consists of the union of sets of functions related to all balls in $\mathbb{Q}_p$.

The set of functions related to a ball has the following form: any of the functions takes values on maximal subballs of the ball (support of the function) equal to (w.r.t. normalization) the $p$-roots of one $\exp(2\pi i m/p)$, $m=0,1,\dots,p-1$, and the values of any function from the set under consideration at the different maximal subballs are different.

The bound of this frame is equal to $p!/(p-1)$.}
\end{lemma}

\noindent{\it Proof}\quad Ball--morphisms which conserve a ball acts by substitutions on the set of maximal subballs of the ball. Applying this observation to the ball $\mathbb{Z}_p$ and the wavelet $\psi$ we get the set of $p!$ functions related to the ball.

Since the parabolic group is transitive on the set of balls in
$\mathbb{Q}_p$ we get that the orbit of the wavelet $\psi$ is the set of functions described above.

Let us show that the constructed set is a tight frame and compute the bound of this frame. To prove that the orbit of the wavelet $\psi$ is a tight frame with the bound $A$ it is sufficient to prove that for the characteristic function $\Omega(|x|_p)$ of the unit ball and the wavelets $\psi^{(N)}$ from the mentioned orbit we will have
$$
\sum_{N}|\langle\Omega(|\cdot|_p),\psi^{(N)}\rangle|^2=A.
$$
The sum is taken over all elements of the orbit.

It is easy to see that the non zero scalar products in the above series have the form
$$
|\langle\Omega(|\cdot|_p),\psi^{(N)}\rangle|= p^{-{\gamma\over 2}},
$$
where the function $\psi^{(N)}$ is related to the ball $J$ with the diameter $p^{\gamma}$,
$\gamma\ge 1$, $J\supset \mathbb{Z}_p$.

Thus
$$
\sum_{N}|\langle\Omega(|\cdot|_p),\psi^{(N)}\rangle|^2
=p!\sum_{\gamma=1}^{\infty}p^{-\gamma}={p!\over p-1}.
$$

The homogeneity of the frame follows from the unitarity of the representation of the parabolic group. $\square$

\section{${\rm mod}\,p$-Affine isometries in many dimensions}

In the multidimensional case one can apply to investigation of isometries the construction of Section 2. In the present Section we formulate the alternative approach --- we define a subgroup of the group of isometries which contains ${\rm mod}\, p$-affine maps.

Let us consider the map $\phi:\mathbb{Q}_p^d\to\mathbb{Q}_p^d$, which is a one--to--one isometry, i.e. an automorphism of the tree of balls in the ultrametric space $\mathbb{Q}_p^d$ belonging to the orispheric group.

For a ball $B\subset \mathbb{Q}_p^d$ we put in correspondence to maximal subballs in $B$ the elements of the module (over the field of residues ${\rm mod}\,p$)  $\mathbb{F}^d_p$. The elements of this module can be considered as vectors $(b_1,\dots,b_d)$, $b_l=0,\dots, p-1$ with the natural structure of the module ${\rm mod}\,p$. We call this module {\it the tangent module} for the ball $B$.

There exists a natural map $T_B: B\to \mathbb{F}^d_p$. Namely for some point $x_B\in B$ we consider the translation of $B$ to itself followed by the dilation and factorization ${\rm mod}\,p\mathbb{Z}^d_p$ of the following form
$$
T_B: B\to \mathbb{F}^d_p,
$$
$$
T_B: x\mapsto (x-x_{B}){\rm diam}(B)\left({\rm mod}\,p\mathbb{Z}^d_p\right).
$$

For the isometry $\phi$ which maps the ball $B$ to the ball $\phi(B)$ we define {\it the tangent map}
$$
\phi_B:\mathbb{F}^d_p\to \mathbb{F}^d_p,
$$
$$
\phi_B=T_{\phi(B)}\circ\phi\circ T_B^{-1}.
$$

Since an isometry maps balls to balls the map $\phi_B$ is defined correctly and is a one--to--one map. This map is defined non uniquely (it depends on $x_B$, $x_{\phi(B)}$).

We shall consider particular isometries which satisfy the following condition.

\begin{definition}{\sl
The isometry $\phi:\mathbb{Q}_p^d\to\mathbb{Q}_p^d$ is called ${\rm mod}\, p$-affine if
for any ball $B\in{\cal T}(\mathbb{Q}_p)$ the corresponding tangent map $\phi_B$ is an affine map of modules over the field of residues ${\rm mod}\, p$, i.e. this map is a combination of  a non-degenerate linear ${\rm mod}\, p$ map and a translation from the module.}
\end{definition}

We denote $A_{p,d}$ the group of ${\rm mod}\, p$-affine isometries of $\mathbb{Q}_p^d$.

Due to surjectivity of $\phi$ the map $\phi_B$ will be non degenerate. Taking different $x_B$ and $x_{\phi(B)}$ we obtain maps $\phi_B$ which differ by translations.

The set of all possible maps $\phi_B$ of the above form constitutes a semidirect product of the factorgroup of the group $O_d$ of norm conserving linear maps of $\mathbb{Q}_p^d\to\mathbb{Q}_p^d$ over the subgroup of matrices from the set $1+{\rm Mat}(p\mathbb{Z}_p^d)$ (where ${\rm Mat}(p\mathbb{Z}_p^d)$ is the set of $d\times d$ square matrices with the matrix elements in $p\mathbb{Z}_p$) and the group of translations $\mathbb{F}^d_p$.
This and the results of \cite{framesdimd} imply the following lemma.

\begin{lemma}{\sl
For the ${\rm mod}\, p$-affine isometry $\phi$ and the ball $B$ with the characteristic function $\Omega(p^{\gamma}\cdot-n)$, which maps to the ball $\phi(B)$ with the characteristic function $\Omega(p^{\gamma}\cdot-n')$, $\gamma\in\mathbb{Z}$, $n,n'\in \mathbb{Q}_p^d/\mathbb{Z}_p^d$, the wavelet $\psi_{J\gamma n}$ supported on the ball $B$ maps to the product of the wavelet $\psi_{\phi_B J,\gamma n'}$ supported on the ball $B'$ and some $p$-root of the unit.}
\end{lemma}

Therefore wavelet bases allow us to study representations of the introduced subgroup of the group isometries of $\mathbb{Q}_p^d$ (equivalently of the orispheric group of automorphisms of the corresponding tree of balls).

\section{Subsets, bases and isometries}

Let us consider the set of linear independent vectors $\{k_1,\dots,k_d\}$, $k_l\in\mathbb{Q}_p^d$ (a basis in $\mathbb{Q}_p^d$). The following lemma describes some sets which can be constructed taking linear combinations of vectors from the basis.

\begin{lemma}\label{setS}{\sl
1) Let the points $z=(z_1,\dots,z_d)$ belong to the ball $|z|_p\le p^{\gamma}$, all vectors in the basis $\{k_1,\dots,k_d\}$ have the norm one and $\{k_1\,{\rm mod}\,p,\dots,k_d\,{\rm mod}\,p\}$ generate $\mathbb{F}^d_p$.
Then there exists a one to one correspondence between the $d$-dimensional ball $z:|z-x_0|_p\le p^{\gamma}$, $x_0\in\mathbb{Q}_p^d$ and the set of points of the form
$$
x_0+\sum_{l=1}^{d}z_lk_l.
$$

2) Let the points $z=(z_1,\dots,z_d)$ belong to the sphere $|z|_p=p^{\gamma}$, all vectors in the basis $\{k_1,\dots,k_d\}$ have the norm one and $\{k_1\,{\rm mod}\,p,\dots,k_d\,{\rm mod}\,p\}$ generate $\mathbb{F}^d_p$.
Then there exists a one to one correspondence between the $d$-dimensional sphere $z:|z-x_0|_p=p^{\gamma}$, $x_0\in\mathbb{Q}_p^d$ and the set of points of the form
$$
x_0+\sum_{l=1}^{d}z_lk_l.
$$

3) Let the points $z=(z_1,\dots,z_d)$ belong to the product of the one dimensional sphere $|z_1|_p=p^{\gamma}$ and the $d-1$-dimensional ball $|z_l|_p\le p^{\gamma}$, $l=2,\dots,d$. Let also the vectors from the basis $\{k_1,\dots,k_d\}$ have the properties: $|k_1|_p=1$, $|k_l|_p=p^{-1}$, $l=2,\dots,d$, and the set of vectors
$$
\{k_1\,{\rm mod}\,p,(p^{-1}k_2)\,{\rm mod}\,p,\dots,(p^{-1}k_d)\,{\rm mod}\,p\}
$$
generates $\mathbb{F}^d_p$.

Then the subset $S\subset\mathbb{Q}_p^d$ of the form
$$
x_0+\sum_{l=1}^{d}z_lk_l
$$
with $x_0\in \mathbb{Q}_p^d$, is the union of $p-1$ balls with the diameter $p^{\gamma-1}$.

Each of the balls (enumerated by $j=1,\dots,p-1$) can be parametrized as
$$
x_0+jp^{\gamma}k_1+ \sum_{l=1}^{d}z_lk_l,\quad |z_1|_p\le p^{\gamma-1},\quad |z_l|_p\le p^{\gamma},\quad l=2,\dots,d.
$$
}
\end{lemma}

\noindent{\it Proof}\quad The unit ball in $\mathbb{Q}_p^d$  with the center in zero is a $\mathbb{Z}_p$-module of the rank $d$.
The set $\{k_1,\dots,k_d\}$ of vectors in $\mathbb{Q}_p^d$ generate the unit ball (as a module) if and only if
the set $\{k_1,\dots,k_d\}$ contains the norm one linear independent vectors and $\{k_1\,{\rm mod}\,p,\dots,k_d\,{\rm mod}\,p\}$ generates $\mathbb{F}^d_p$.

An arbitrary ball is a translation and dilation of the unit ball with the center in zero. This proves the first statement of the lemma.

A sphere is a difference of two balls. This implies the second statement of the lemma.

Statement 3 of the lemma: consider the set $\{jp^{\gamma}\}$, $j=1,\dots,p-1$ of representatives in the maximal subballs in the sphere $|z_1|_p=p^{\gamma}$. Then we can apply the first statement of the lemma to any of the balls with the centers in $x_0+jp^{\gamma}k_1$, the basis
$\{k_1,p^{-1}k_2\dots,p^{-1}k_d\}$ and the set $|z|_p\le p^{\gamma-1}$.

This finishes the proof of the lemma. $\square$

\begin{remark}\label{reconstruct_basis}{\rm Equivalently, if we have the set $S$ (which is a subset of the ball $B={\rm sup}(S)$ -- of the minimal ball which contains $S$) described in the part 3 of the above lemma, we can recover the basis $\{k_1,\dots,k_d\}$ in the following way. Taking two different maximal subballs in $S$ and choosing some points $y_0$, $y_1$ in these balls, we put $k_1={\rm diam}(B)(y_1-y_0)$ (with this choice $|k_1|_p=1$). Then we choose arbitrary vectors $k_2$, \dots, $k_d$ in such a way that $|k_l|_p=p^{-1}$, $l=2,\dots,d$, $\{k_1,p^{-1}k_2\dots,p^{-1}k_d\}$ is a basis in $\mathbb{Q}_p^d$ and the set $\{k_1\,{\rm mod}\,p,(p^{-1}k_2)\,{\rm mod}\,p,\dots,(p^{-1}k_d)\,{\rm mod}\,p\}$ generates $\mathbb{F}^d_p$.

We also recover the point $x_0\in B$ as follows: taking the obtained vector $k_1$ and any point $y_0\in S$ we consider the set of $p$ points in $B$ of the form $y_j=y_0+j\,{\rm diam}^{-1}(B)k_1$, $j=0,\dots,p-1$. Then $p-1$ of the described points will lie in $S$ and one of the points will not belong to $S$. We denote this point $y_A$. Then taking the maximal subball in $B$ containing $y_A$ we choose an arbitrary point $x_0$ in this subball.

The basis $\{k_1,\dots,k_d\}$ and $x_0$ so constructed are not unique but if we apply to any of the constructed bases and points $x_0$ the procedure of the part 3 of lemma \ref{setS} we reproduce the set $S$ unambiguously.
}
\end{remark}

\begin{remark}\label{basis_in_tangent_space}{\rm Using the previous remark one can consider the following one--to--one correspondence between the sets $S$ described in the Statement 3 of lemma \ref{setS} and pairs $(\mathbb{F}_p^{*}k_1,B_0)$ where $k_1\in\mathbb{F}^d_p$, $k_1\ne 0$, $\mathbb{F}_p^{*}$ is the set of non zero elements in the field $\mathbb{F}_p$ of residues ${\rm mod}\, p$ and $B_0\in\mathbb{Q}_p^d/p\mathbb{Z}_p^d$ (equivalently $B_0$ can be considered as the maximal subball in $B=\,{\rm sup}\,(S)$ which contains the point $x_0$).

To reproduce the set $S$ we  consider $k_1$ as a vector in $\mathbb{Q}_p^d$ (taking an arbitrary representative in the equivalence class), then we take an arbitrary set $\{k_2,\dots,k_d\}$ such that $k_1\bigcup\{p^{-1}k_2,\dots,p^{-1}k_d\}$ generate $\mathbb{Z}_p^d$ as a $\mathbb{Z}_p$-module. Then we choose an arbitrary point $x_0\in B_0$.

The obtained basis $\{k_1,\dots,k_d\}$ together with the point $x_0$ generate the set $S$ with the help of the part 3 of lemma \ref{setS}. The set $S$ obtained in this way  will depend only on $\mathbb{F}_p^{*}k_1$, $k_1\ne 0\in\mathbb{F}^d_p$ and the maximal subball $B_0$ in $B$,   $x_0\in B_0$.

The pair $(k_1,B_0)$ can be considered as a vector in the affine space $\mathbb{F}^d_p$ (i.e. a vector together with the point to which the vector is attached).
}
\end{remark}

The statement of the next lemma  describes the images of the sets $S$ described in the part 3 of lemma \ref{setS} with respect to a ${\rm mod}\, p$-affine isometry.

\begin{lemma}\label{tangent_map}{\sl
Let $\phi$ be a ${\rm mod}\, p$-affine isometry and $S$ be a subset of the ball $B\subset\mathbb{Q}_p^d$ corresponding to a non-zero vector in the affine space $(k_1,B_0)\in \mathbb{F}^d_p$. Then the subset $\phi(S)$ of the ball $\phi(B)$ will correspond to the vector $\phi(k_1,B_0)=(\phi(k_1),\phi(B_0))$ in the affine space.

Equivalently, for the set $S$ corresponding to the point $x_0\in B$ and the basis $\{k_1,\dots,k_d\}$ in $\mathbb{Q}_p^d$, $\{k_1\,{\rm mod}\,p,(p^{-1}k_2)\,{\rm mod}\,p,\dots,(p^{-1}k_d)\,{\rm mod}\,p\}$ generates $\mathbb{F}^d_p$, and the ${\rm mod}\, p$-affine isometry $\phi$, the set $\phi(S)$ will be a set belonging to the same family of sets (described in the part 3 of lemma \ref{setS}) and will correspond to the point $\phi(x_0)$ and the basis $\{\phi(k_1),\dots,\phi(k_d)\}$.}
\end{lemma}

\noindent{\it Proof}\quad The set $S$ is a union of $p-1$ balls. Since $\phi$ is an isometry it maps disjoint balls to disjoint balls with the same diameters. Therefore it is sufficient to check the action of $\phi$ on the centers of the mentioned balls. The centers of the balls in $S$ are the points $x_0+jp^{\gamma}k_1$, $j=1,\dots,p-1$. The images of these points by ${\rm mod}\, p$-affinity of $\phi$ are
$$
\phi(x_0+jp^{\gamma}k_1)=\left(\phi(x_0)+jp^{\gamma}\phi_{B}(k_1)\right)\,{\rm mod}\,p^{\gamma+1}.
$$
Here we identify $k_1\in\mathbb{Q}_p^d$, $|k_1|_p=1$ with its image in $\mathbb{F}^d_p$ (this makes sense since we consider the above identity $\,{\rm mod}\,p^{\gamma+1}$).

This set has the form of a set described in the part 3 of lemma \ref{setS}.
This finishes the proof of the lemma. $\square$

\begin{remark}
{\rm The above lemma can be understood as follows. We have the action of the group $A_{p,d}$ of ${\rm mod}\, p$-affine isometries in $\mathbb{Q}_p^d$. This action generates the action of $A_{p,d}$  by automorphisms of the tree ${\cal T}(\mathbb{Q}_p^d)$.

Each vertex $B$ in ${\cal T}(\mathbb{Q}_p^d)$ is equipped with the module $\mathbb{F}^d_p$ (the elements of the module are in one to one correspondence to maximal subballs in $B$).  We consider this family of modules as a fiber bundle over ${\cal T}(\mathbb{Q}_p^d)$.

To each $B\in {\cal T}(\mathbb{Q}_p^d)$ we associate the corresponding set $S$ (equivalently, family of vectors $(\mathbb{F}_p^{*}k_1,B_0)$ where $(k_1,B_0)$ is a non-zero vector in the affine space $\mathbb{F}^d_p$). The family of sets $S$ corresponding to all vertices in the tree ${\cal T}(\mathbb{Q}_p^d)$ can be considered as an equivalence class of affine vector fields on ${\cal T}(\mathbb{Q}_p^d)$ which do not take zero values where two affine vector fields are equivalent if one can be made equal to the other by multiplication by a $\mathbb{F}_p^{*}$-valued function.

Lemma \ref{tangent_map} states that the tangent map for a transformation from $A_{p,d}$ maps the fiber bundle on ${\cal T}(\mathbb{Q}_p^d)$ of affine vector fields with values in $\mathbb{F}^d_p$ to itself.

}
\end{remark}

\begin{remark}
{\rm The vectors $k_1\ne 0\in\mathbb{F}^d_p$ are in one to one correspondence to maximal subballs in the unit sphere in $\mathbb{Q}_p^d$. This set of balls is parametrized by  $J$ in (\ref{psi_J}), (\ref{J}). This is exactly the parameter on the set of wavelets associated to a ball in $\mathbb{Q}_p^d$.

The correspondence is straightforward --- we put $J=k_1$ where $k_1\in\mathbb{F}^d_p$ is considered as a vector in $\mathbb{Q}^d_p$ with the coordinates equal to the corresponding residues in the finite field of residues.
Therefore the mentioned above fiber bundle can be considered as a fiber bundle of wavelets.}
\end{remark}

\section{Operators and isometries}

Let us denote by $B(x_0,p^{\gamma})$ the ball with the center in $x_0$ and the diameter $p^{\gamma}$. Consider the vector field  ${\bf k}=\{(k_1(B(x_0,p^{\gamma}))\}$ on the tree of balls ${\cal T}(\mathbb{Q}_p^d)$. Here $k_1$ is a non-zero vector in $\mathbb{F}^d_p$.

Correspondingly (by remark \ref{basis_in_tangent_space}) ${\bf k}$ describes the family of bases described in the part 3 of lemma \ref{setS}.

Using the above notations let us introduce the integral operator of the following form.

\begin{definition}{\sl
The following operator will be called a pseudodifferential vector field
\begin{equation}\label{vector_field}
D_{F(\cdot),\bf k(\cdot)}f(x)=\sum_{\gamma\in\mathbb{Z}}\int_{|z_1|_p=p^{\gamma}\atop |z_l|\le p^{\gamma},l=2,\dots,d} F(B(x,p^{\gamma}))\left[f(x)-f\left(x+\sum_{l=1}^{d}z_lk_l(B(x,p^{\gamma})\right)\right]dz_1\dots dz_d.
\end{equation}
Here $F(B(x,p^{\gamma}))$ is a complex valued integration kernel defined on the tree of balls ${\cal T}(\mathbb{Q}_p^d)$.}
\end{definition}

This operator can be considered as an operator of pseudodifferentiation along the vector field ${\bf k}$ on the tree of balls ${\cal T}(\mathbb{Q}_p^d)$.

\begin{remark}{\rm For fixed $x$ and $\gamma$ the integration in (\ref{vector_field}) runs over the set $S\subset B(x,p^{\gamma})$ which corresponds to the point $x$ and to the basis $\{k_1,\dots,k_d\}$. Sets of integration of this kind are in one--to--one correspondence to vectors $(k_1,B_0)$ in the affine space  $\mathbb{F}^d_p$, where $B_0$ is the maximal subball in $B(x,p^{\gamma})$ containing the point $x$ and $k_1\ne 0\in\mathbb{F}^d_p$ corresponds to the basis $\{k_1,\dots,k_d\}$.
}
\end{remark}

\begin{remark}
{\rm The above formula (\ref{vector_field}) may look complicated, but this is a direct generalization of a standard pseudodifferential operator (which can be diagonalized by the Fourier transform). For a standard pseudodifferential operator we have the expression
\begin{equation}\label{standard}
Df(x)=\int_{\mathbb{Q}_p^d} F(|x-y|_p)(f(x)-f(y))dy=
$$
$$
=\sum_{\gamma\in\mathbb{Z}}\int_{|z|_p=p^{\gamma}}F(p^{\gamma})\left[f(x)-f\left(x+\sum_{l=1}^{d}z_lk_l(B(x,p^{\gamma})\right)\right]dz_1\dots dz_d
\end{equation}
where the basis $\{k_l(B)\}=\{k_1,\dots,k_d\}$ for any ball $B$ contains vectors having norm one.

The main differences between the above standard definition and formula (\ref{vector_field}) are the following:

1) The integration domain in (\ref{vector_field}) (for a fixed $\gamma$) is the product of a one dimensional sphere and a ball, instead of a $d$-dimensional sphere in
(\ref{standard});

2) The basis $\{k_1,\dots,k_d\}$ in (\ref{vector_field}) contains one vector of norm one and the other vectors are of norm $p^{-1}$, instead
of basis containing vectors all of norm one as in (\ref{standard}) (moreover operators (\ref{standard}) do not depend on the choice of the corresponding bases);

3) In (\ref{standard}) we might consider the integration kernel $F(B(x,p^{\gamma}))$ as in (\ref{vector_field}). Actually a family of operators of the form (\ref{standard}) with kernels $F(B(x,p^{\gamma}))$ was considered in \cite{nhoper} for the $p$-adic case and in \cite{Izv}, \cite{ACHA}, \cite{MathSbornik} for general locally compact ultrametric spaces.
}
\end{remark}

For a ${\rm mod}\, p$-affine isometry $\phi$ we denote by $\phi({\bf k})$ the set of bases $\{\phi(k_1),\dots,\phi(k_d)\}$ corresponding to balls $B$.

The next theorem describes the transformation with respect to a ${\rm mod}\, p$-affine isometry of the operator $D_{F(\cdot),{\bf k}(\cdot)}$ given by (\ref{vector_field}) corresponding to the integration kernel $F(B)$ (complex valued function on ${\cal T}(\mathbb{Q}_p^d)$) and the vector field ${\bf k}(B)$ on ${\cal T}(\mathbb{Q}_p^d)$ with non-zero values in $\mathbb{F}^d_p$.

\begin{theorem}{\sl
Let $\phi$ be a ${\rm mod}\, p$-affine isometry. Then
\begin{equation}\label{change_variable}
\Phi\circ D_{F(\cdot),\bf k(\cdot)} f(x)=D_{F(\phi(\cdot)),\bf \phi^{-1}k(\phi(\cdot))}\circ \Phi f(x).
\end{equation}
}
\end{theorem}

Here we denote $\Phi$ the map on functions on $\mathbb{Q}_p^d$ which corresponds to the  ${\rm mod}\, p$-affine isometry $\phi$:
$$
\Phi[f](x)=f(\phi(x)).
$$
Formula (\ref{change_variable}) is the analogue of the transformation of a vector field considered as a pseudodifferential operator with respect to a change of coordinates.

\bigskip

\noindent{\it Proof}\quad Let us consider
\begin{equation}\label{proof}
D_{F(\cdot),\bf k(\cdot)}\circ \Phi f(x)=D_{F(\cdot),\bf k(\cdot)}f(\phi(x))=
$$
$$
=\sum_{\gamma\in\mathbb{Z}}\int_{|z_1|_p=p^{\gamma}\atop |z_l|\le p^{\gamma},l=2,\dots,d} F(B(x,p^{\gamma}))\left[f(\phi(x))-f\left(\phi\left(x+\sum_{l=1}^{d}z_lk_l(B(x,p^{\gamma})\right)\right)\right]dz_1\dots dz_d.
\end{equation}

Consider the map
$$
\phi:S\to\phi(S),
$$
where the set $S$ is described by the part 3 of lemma \ref{setS}:
$$
S=\{x+\sum_{l=1}^{d}z_lk_l(B(x,p^{\gamma}), \quad |z_1|_p=p^{\gamma},|z_l|\le p^{\gamma},l=2,\dots,d\}
$$

By lemma \ref{tangent_map} the set $\phi(S)$ has the above form where instead of $x$ and $\{k_1,\dots,k_d\}$ we have to use  $\phi(x)$ and $\{\phi(k_1),\dots,\phi(k_d)\}$

$$
\phi\left(x+\sum_{l=1}^{d}z_lk_l(B(x,p^{\gamma})\right)=\phi(x)+\sum_{l=1}^{d}z'_l\phi(k_l(B(x,p^{\gamma})),
$$
where $z'_l$ belong to the same subset of $\mathbb{Q}_p^d$.

Moreover since the map $\phi$ is an isometry it conserves the measure, therefore
$$
dz'_1\dots dz'_d=dz_1\dots dz_d.
$$

We get for (\ref{proof}) the expression
$$
D_{F(\cdot),\bf k(\cdot)}\circ \Phi f(x)=
$$
$$
=\sum_{\gamma\in\mathbb{Z}}\int_{|z_1|_p=p^{\gamma}\atop |z_l|\le p^{\gamma},l=2,\dots,d} F(B(x,p^{\gamma}))\left[f(\phi(x))-f\left(\phi(x)+\sum_{l=1}^{d}z_l\phi(k_l(B(x,p^{\gamma}))\right)\right]dz_1\dots dz_d.
$$

Analogously
$$
\Phi\circ D_{F(\cdot),\bf k(\cdot)} f(x)=
$$
$$
=\sum_{\gamma\in\mathbb{Z}}\int_{|z_1|_p=p^{\gamma}\atop |z_l|\le p^{\gamma},l=2,\dots,d} F(B(\phi(x),p^{\gamma}))\left[f(\phi(x))-f\left(\phi(x)+\sum_{l=1}^{d}z_lk_l(B(\phi(x),p^{\gamma})\right)\right]dz_1\dots dz_d.
$$

We get
$$
\Phi\circ D_{F(\cdot),\bf k(\cdot)} f(x)=D_{F(\phi(\cdot)),\bf \phi^{-1}k(\phi(\cdot))}\circ \Phi f(x).
$$

This finishes the proof of the theorem.
$\square$

\section{Appendix}

In this section we put some definitions and notations used in $p$-adic analysis.

We denote $\mathbb{F}_p$ the field  of residues ${\rm mod}\, p$ and $\mathbb{F}_p^{*}$ the multiplicative group of non zero elements in $\mathbb{F}_p$.

For a $p$-adic number $x$ with the expansion over the degrees of $p$ (where $p$ is prime) of the form
\begin{equation}\label{x}
x=\sum_{j=\gamma}^{\infty}x_{j}p^{j},\quad x_j=0,\dots,p-1,\quad \gamma\in\mathbb{Z},\quad x_{\gamma}\ne 0
\end{equation}
its $p$-adic norm is $|x|_p=p^{-\gamma}$.

The norm in $\mathbb{Q}_p^{d}$ is introduced as follows: for $x=(x_1,\dots,x_d)\in\mathbb{Q}_p^{d}$ the norm is defined as
$$
|x|_p={\rm max}_{l=1,\dots, d} |x_l|_p.
$$

$p$-Adic wavelets were introduced in \cite{wavelets}. Using the notations of \cite{framesdimd} the (multidimensional) wavelet basis $\{\psi_{\gamma n J}\}$ in the space $L^2(\mathbb{Q}_p^{d})$ has the form of the set of functions of the form

$$
\psi_{\gamma n J}(x)=p^{-{d\gamma\over
2}}\psi_{J}(p^{\gamma}x-n),\qquad x\in \mathbb{Q}_p^d,\quad
\gamma\in\mathbb{Z},\quad n\in \mathbb{Q}_p^d/\mathbb{Z}_p^d,
$$
$$
n=\left(n^{(1)},\dots,n^{(d)}\right),\qquad
n^{(l)}=\sum_{i=\beta_l}^{-1}n^{(l)}_{i}p^{i},\quad
n_i^{(l)}=0,\dots,p-1,\quad \beta_l\in\mathbb{Z}_{-}.
$$
where
\begin{equation}\label{psi_J}
\psi_J(x)=\chi(p^{-1}Jx)\Omega(|x|_p),\qquad x,k\in \mathbb{Q}^d_p,\quad Jx=\sum_{l=1}^d J_l x_l,
\end{equation}
Here
\begin{equation}\label{J}
J=\left(j_1,\dots,j_d\right),\qquad j_l=0,\dots,p-1,
\end{equation}
where at least one of $j_l$ is not equal to zero. The function $\Omega(|x|_p)$ is the characteristic function of the unit ball and $\chi$ is the complex valued character of a $p$-adic argument which for $x$ with the expansion (\ref{x}) takes the form
$$
\chi(x)=\exp\left(2\pi i \sum_{j=\gamma}^{-1}x_{j}p^{j}\right).
$$

The Vladimirov operator of $p$-adic fractional differentiation is the pseudodifferential operator of the form
$$
D^{\alpha}f(x)=  F^{-1}[|k|_p^{\alpha} F[f]](x),
$$
where $F$ is the $p$-adic Fourier transform
$$
F[f](k)=\int \chi(kx)f(x)dx
$$
and $F^{-1}$ is the inverse Fourier transform, the integration is taken with respect to the $p$-adic Haar measure.

For $\alpha>0$ this operator can be put into the form
$$
D^{\alpha} f(x)=\frac{1}{\Gamma_p(-\alpha)}
\int_{Q_p}\frac{f(x)-f(y)}{|x-y|_p^{1+\alpha}}dy,\quad \Gamma_p(-\alpha)={p^{\alpha}-1\over 1-p^{-1-\alpha}}.
$$

Ultrametric spaces are dual to trees with some partial order. Below we describe
some part of the duality construction.

For a (complete locally compact) ultrametric space $X$ we consider the set ${\cal T}(X)$,
which contains all the balls in $X$ of nonzero diameters, and the
balls of zero diameter which are maximal subbals in balls of nonzero
diameters. This set possesses a natural structure of a partially ordered
tree. The partial order in ${\cal T}(X)$ is defined by inclusion of balls.

Two vertices $I$ and $J$ in ${\cal T}(X)$ are connected by an
edge if the corresponding balls are ordered by inclusion, say
$I\supset J$ (i.e. one of the balls contains the other), and there
are no intermediate balls between $I$ and $J$.

On the tree ${\cal T}(X)$ we have the natural increasing
positive function which puts in correspondence to any vertex the
diameter of the corresponding ball.

Assume now that we have a partially ordered tree ${\cal T}$, satisfying the conditions:

1) Graph ${\cal T}$ is a tree, i.e. for any pair of vertices there exists a finite path in ${\cal T}$ which connects these vertices and ${\cal T}$ does not contain cycles.

2) Each vertex in ${\cal T}$ is incident to a finite set of edges.

3) For any finite path in ${\cal T}$ there exists a unique maximal vertex in this path.

Let us choose an arbitrary  positive
increasing (w.r.t. the partial order) function $F$ on this tree. Then we define the ultrametric
on the set of vertices of the tree ${\cal T}$ as follows:
$$
d(I,J)=F({\rm sup}(I,J)),
$$
where ${\rm sup}(I,J)$ is the supremum of vertices $I$,
$J$ with respect to the partial order. The vertex ${\rm sup}(I,J)$ coincides with the above mentioned unique maximal vertex in the path $IJ$.

Then we take a completion of the set of vertices with respect to the
defined ultrametric and eliminate from the completion all the inner
points of the tree (a vertex of the tree is inner if it does not
belong to the border of the tree). We denote the obtained space
$X({\cal T})$, this space is ultrametric and locally compact.

\begin{definition} {\sl
The set of vectors $\{f_n\}$ in the Hilbert space ${\cal H}$ is a
frame, if there exist positive constants $A,B>0$, such that for each
vector $g\in {\cal H}$ the following inequality is satisfied:
$$
A\|g\|^2\le\sum_{n} |\langle g,f_n\rangle|^2 \le B\|g\|^2.
$$
}
\end{definition}

The constants $A$ and $B$ are called the lower and the upper bounds
of the frame correspondingly. A frame is tight if the frame bounds
$A$ and $B$ are equal. A frame is uniform if all elements have equal
norms.

\bigskip

\noindent{\bf Acknowledgments}\qquad One of the authors (S.K.) would
like to thank I.V.Volovich, V.S.Vla\-di\-mi\-rov and E.I.Zelenov
for fruitful discussions and valuable comments. He gratefully
acknowledges being partially supported by the grant DFG Project 436
RUS 113/809/0-1 and DFG Project 436 RUS 113/951, by the grants of
the Russian Foundation for Basic Research
RFFI 08-01-00727-a and RFFI 09-01-12161-ofi-m, by the grant of the President of Russian
Federation for the support of scientific schools NSh 3224.2008.1 and
by the Program of the Department of Mathematics of the Russian
Academy of Science ''Modern problems of theoretical mathematics'',
and by the program of Ministry of Education and Science of Russia
''Development of the scientific potential of High School, years of
2009--2010'', project 3341. He is also grateful to IZKS (the
Interdisciplinary Center for Complex Systems) of the University of
Bonn for kind hospitality.

\end{document}